\documentclass[12pt]{article}

\usepackage[all]{xy}

\usepackage[mathscr]{eucal}

\usepackage{amssymb}
\usepackage{amsthm}
\usepackage{amsmath}

\newcommand{\R}{\mathbb{R}}
\newcommand{\C}{\mathbb{C}}

\def\cC{{{\Cal C}}}
\def\cK{{{\cal K}}}

\def\cT{{{\cal{T}}}}
\def\cC{{{\cal C}}}

\def\IR{{{\Bbb R}}}

\def\IZ{{{\Bbb Z}}}

\def\Hom{{\text{Hom}}}

\def\C{{\text C}}
\def\rM{{\text M}}
\def\K{{\frak K}}
\def\smallsim{^{\frak{\sim}}}
\renewcommand{\thefootnote}

\begin{document}


\title{A remark on invariants for C*-algebras
of stable rank one}

\author{Alin  Ciuperca and George A. Elliott}
\maketitle

\begin{abstract}

It is shown that, for a C*-algebra of stable rank one (i.e., in which the invertible elements are dense), two
well-known isomorphism invariants, the Cuntz semigroup and the Thomsen semigroup, contain the same information.
More precisely, these two invariants, viewed appropriately, determine each other in a natural way.
\end{abstract}

\bigskip
{\bf 1.}\ \ In [5], Cuntz introduced an invariant for C*-algebras---very much analogous to the purely algebraic
invariant K$_0$---which has come to be known as the Cuntz semigroup, and which has recently received
considerable attention (see [21], [2], [4], and [12]). (In [21], Toms showed that the Cuntz semigroup is of
crucial importance in distinguishing C*-algebras.)

In [20], Thomsen introduced a related (somewhat more complicated) invariant, which we shall refer to as the
Thomsen semigroup. Thomsen showed that his invariant was complete in the case of separable approximate interval
algebras---sometimes called AI algebras.

The purpose of this paper is to point out certain relations between these two invariants, in the case of stable
rank one, which follow from the results of [4] (see Theorems 4 and 10). (Roughly speaking, in this case the two
invariants coincide.) In particular, it will follow from our remarks that the Cuntz semigroup is a complete
invariant for AI algebras (see Theorem 11). \vspace{8pt}

\footnote{The research of the second author was supported by a grant from the Natural Sciences and Engineering Research Council of Canada. This work was completed during the 2007 Thematic Program on Operator Algebras and Applications at Fields Institute for Research in Mathematical Science.}

{\bf 2.}\ \ Before making a precise comparison of the Cuntz and Thomsen semigroups it is of course necessary to
review the definitions of these invariants. In particular, it is important to recall the mathematical
structures---in more modern language, the categories---involved, as these objects are both much more than just
abstract semigroups. The elements of the Thomsen semigroup are easily described; they are just the approximate
unitary equivalence classes of C*-algebra maps (*-homomorphisms) from the reference C*-algebra $\C_0 (]0, 1])$
to the given C*-algebra under study---say $A$---or to the stabilization of $A$, the tensor product $A\otimes \K$
where $\K=\K (l^2)$ is the C*-algebra of compact operators, if $A$ is not stable. Recall that two maps from a
C*-algebra $B$ into $A$ are said to be approximately unitarily equivalent if they are arbitrarily close modulo
inner automorphisms of $A$ on any finite subset of $B$, with respect to the norm of $A$ (where by an inner
automorphism of $A$ we mean an automorphism determined by a unitary element of the *-algebra obtained by
adjoining a unit to $A$)---in other words, if the closures of their orbits under unitary equivalence, with
respect to the topology of pointwise convergence in norm, are equal. This topology gives rise to a topology on
the quotient space---the space of closures of unitary orbits. In the case $B=\C_0(]0, 1])$ this topology arises
in a natural way from a complete metric, which passes to a complete metric on the quotient.

As far as the semigroup structure of this quotient space is concerned, the simplest definition is to choose  an
isomorphism of $\K$ with $\K\otimes \rM_2$, the C*-algebra of $2\times 2$ matrices over $\K$, and for two given
C*-algebra maps $\varphi$ and $\psi$ from $B$ to $A\otimes \K$ to consider the block diagonal map $\varphi
\oplus \psi$ from $B$ to $A\otimes \K \otimes \rM_2$, i.e., the map
\[ b\mapsto \left(\begin{array}{cc}
\varphi (b)  &0 \\
0 & \psi (b)\end{array}\right),
\]
and to consider this as a map from $B$ to $A\otimes \K$, with respect to the chosen isomorphism of $\K\otimes
\rM_2$ with $\K$. The closure of the unitary orbit of the sum of $\varphi$ and $\psi$ defined in this way
depends only the closures of the unitary orbits of $\varphi$ and $\psi$---given that any two isomorphisms of
$\K\otimes \rM_2$ with $\K$ are approximately unitarily equivalent (with respect to unitaries in the algebra
$\K$ with unit adjoined). Commutativity and associativity, at the level of orbit closures, are easily verified,
and so one has the structure of abelian semigroup with zero. Uniform continuity of addition is clear, and so the
overall structure (in the case that $B$ is separable) is that of complete metric abelian semigroup with zero.
(In the case $B=\C_0 (]0, 1])$ one can even say that the metric is natural---in general this applies only to the
uniform structure.)

In addition to the structure of metric abelian semigroup with zero, the Thomsen invariant also admits a natural
(right) action of the (multiplicative) monoid of endomorphisms of the C*-algebra $\C_0(]0, 1])$ (consisting
of---or, rather, arising from---composition on the domain side). This action is compatible with the semigroup
structure (since right multiplication by an endomorphism of $\C_0(]0, 1])$ distributes over addition of two
homomorphisms with orthogonal images from $\C_0(]0, 1])$ to another C*-algebra).

As morphisms in the category $\cT h$ of complete metric abelian semigroups with zero, with a right action of the
endomorphism monoid of $\C_0(]0, 1])$ respecting addition together with the uniform structure, following [19]
let us choose contractive additive maps which respect (i.e.~commute with) the action of the endomorphism monoid.
(Here, by contractive we mean not strictly increasing any distance.) As pointed out in [19], the functor which
to a C*-algebra $A$ associates its Thomsen semigroup, say ${\cT h} A$ in the category $\cT h$, preserves
sequential inductive limits (which always exist in the category $\cT h$).

It follows, as shown in [20], that the functor $\cT h$ exactly describes the approximate unitary equivalence
classes of maps from a separable AI C*-algebra $A$ into an arbitrary C*-algebra $B$. Hence, as follows for
instance from Theorem 3 of [9], in view of Section 4.3 of [9] (and was of course also shown in [20]), the
functor $\cT h$ determines isomorphism of (separable) AI algebras---in other words is a classification functor
for these algebras in the sense of [9].

In [20], the question was raised of investigating the relationship between this classification result and the
classification of simple AI algebras obtained by the second author of the present paper around the same time in
[8], using K$_0$ and traces. The result of the present paper may be regarded as a step in this direction---given
that the Cuntz semigroup for simple AI algebras (among many other simple C*-algebras!) has been calculated in
[2]. (A computation of the Thomsen invariant for a simple AI algebra was also given in [9], using directly the
methods of [8].) Note that a classification of a wide class of non-simple AI algebras has been obtained in [11]
(roughly speaking, the case that K$_1$ of every closed two-sided ideal is zero). (See [19] for an earlier such
result, concerning the rather more special case that every closed two-sided ideal is generated by projections.)
It would be very interesting to calculate the Cuntz semigroup for an arbitrary AI algebra in terms of the
tracial and K$_0$ invariants considered in [11]---as this would extend the classification of [11] to the class
of all AI algebras. \vspace{8pt}

{\bf 3.}\ \ The Cuntz semigroup, defined as an ordered abelian semigroup in [5], was shown recently in [4] to
have additional structure, which will be necessary for our calculation of the Thomsen semigroup (in the case of
a C*-algebra of stable rank one). More precisely, the Cuntz semigroup was shown in [4] to belong to the category
$\cC u$ of ordered abelian semigroups (with zero) with the two additional order-theoretic properties

\begin{itemize}
\item[(i)] the supremum of any increasing sequence (equivalently, of any countable upward directed set) exists,
and

\item[(ii)] each element is the supremum of a rapidly increasing sequence, i.e., an increasing sequence such
that each term is compactly contained in the next, in the sense that if the second element is less than or equal
to the supremum of an increasing sequence then the first element is eventually less than or equal to the terms
of this sequence,
\end{itemize}
\noindent and, furthermore, the property that the operation of passing to the supremum of an increasing sequence
and the relation of compact containment (written $<<$) are compatible with addition (i.e., $\sup (a_n +b_m)=
\sup a_n+\sup b_n$, and $a_i << b_i$ for $i=1, 2$ implies $a_1 +a_2 << b_1 +b_2$).

In order to ensure that the Cuntz semigroup of a C*-algebra $A$, denoted by $\cC u A$ in [4] to emphasize that
it belongs to the category $\cC u$, in fact belongs to this category, it was necessary to revise Cuntz's
definition slightly---as otherwise sometimes not all increasing sequences (only bounded ones) would have
suprema. Namely, at the very beginning (just as above for the Thomsen semigroup), the C*-algebra $A$ was
replaced by $A\otimes \K$. Cuntz's definition then yields the semigroup denoted by $\cC u A$ in [4] (see [4],
Appendix 6), but in order to prove that $\cC u A$ belongs to the category $\cC u$ it was convenient in [4]
(indeed, almost necessary---and at the least natural) to rephrase the definition given in [5]. (And in doing
this it became unnecessary to mention the C*-algebra $\K$!)  Namely (and in analogy with K$_0$!), the semigroup
$\cC u A$ was defined as the set of equivalence classes of countably generated Hilbert $A$-modules, with respect
to a moderately simple equivalence relation---amounting in the case that $A$ has stable rank one (as shown in
Theorem 3 of [4]) simply to isomorphism---and with respect to a moderately simple (pre-) order
relation---amounting again in the case that $A$ has stable rank one simply to inclusion!

Since, for the purposes of the present note, we are interested only in the case that $A$ has stable rank one,
let us just take $\cC u A$ to be the ordered semigroup described above---Hilbert $A$-modules ordered by
inclusion and then considered only up to isomorphism, with addition direct sum (which of course preserves both
inclusion and isomorphism). By Theorems 3 and 2 of [4], $\cC uA$ defined in this way belongs to the category
$\cC u$. By Theorem 2 of [4], the functor $\cC u$ preserves sequential inductive limits (which, as also assured
by this theorem, always exist in the category $\cC u$). (In particular, this makes it in principle simple to
compute the Cuntz semigroup of a (separable) AI algebra, and also simple to establish our main result, the
computation of the Thomsen semigroup in terms of the Cuntz semigroup, for AI algebras.) \vspace{8pt}

{\bf 4.}\ \ By the functorial nature of the Cuntz semigroup, as a functor from C*-algebras to the category $\cC
u$ (see Theorem 2 of [4]), to any C*-algebra map
$$\C_0(]0, 1])\to A\otimes \K$$
there is associated a map
$$\cC u \C_0(]0, 1])\to \cC u A$$
in the category $\cC u$. (This holds without any hypothesis on the C*-algebra $A$.) Furthermore, from the
equivalence of the present definition of $\cC u A$ with that of [5] (see Appendix 6 of [4]), it follows
immediately that two approximately unitarily equivalent maps from $\C_0(]0, 1])$ to $A$ give rise to the same
map from $\cC u \C_0(]0, 1])$ to $\cC u A$. (In the case of exact unitary equivalence this can also be seen
purely algebraically by studying the action of an inner automorphism of $A$ on a Hilbert $A$-module---it
transforms it into an isomorphic one for any $A$.) One thus has a (semigroup) map
$$\cT h A\to \Hom (\cC u \C_0(]0, 1]), \cC u A).$$

In fact, since $\cT h A$ may also be identified with approximate unitary equivalence classes of unital
C*-algebra maps
$$\C([0, 1])\to (A\otimes \K)^{\sim},$$
i.e. from the unitization of $\C_0(]0, 1])$ to the unitization of $A\otimes \K$, one also has a (semigroup) maps
$$\cT h A\to \Hom_1 (\cC u \C([0, 1],\ \cC u(A\otimes \K)^{\sim}),$$
where the subscript one denotes unital maps (i.e., maps respecting the images of the units of the algebras in
the Cuntz semigroups).

One might think that these two semigroups of maps between Cuntz semigroups might be (essentially) exactly the
same, and we shall show that every map in the latter semigroup is indeed determined by its restriction to the
subsemigroup $\cC u \C_0(]0, 1])$ of its domain, but we have been unable to show that every map from $\cC u
\C_0(]0, 1])$ to $\cC u A\otimes \K$ extends to a map between $\C([0, 1])$ and $ \cC u(A\otimes \K)^{\sim}$.
\smallskip

\noindent {\bf Theorem 4.1.}\ \ {\it{If $A$ has stable rank one,  then the map from $\cT h A$ to\\
${\text{Hom$_1$}}(\cC u C([0, 1]), \cC u(A\otimes \K){\smallsim})$ just constructed is an isomorphism of
semigroups, compatible with the natural (right) actions of the monoid of endomorphisms of the C*-algebra $\C([0,
1])$. It is contractive with respect to an intrinsic metric on the codomain semigroup (in other words, depending
on $A$ only through the  ordered semigroup $\cC u A$---in fact only through the order structure of $\cC u A$!).
The inverse of this map is uniformly continuous.}}

\begin{proof}
For simplicity, let us consider the more elementary semigroup of maps $\Hom (\cC u \C_0(]0, 1]), \cC u A)$, and
just bear in mind that at some point we will need to assume that a given map in this semigroup extends.

Let us first show that with respect to a suitable (natural) metric on the codomain semigroup, the map is
contractive but does not decrease distances by more than a factor of eight. As it is clear that the map is
additive and compatible with the monoid action, it then remains to prove surjectivity (and to do this it is
enough to prove that the image is dense).

The metric is simple enough to define, once one recalls the computation of $\cC u \C_0(]0, 1])$ given
(essentially) in [14]---see Theorem 10 below---, namely, as the ordered semigroup of lower semicontinuous
functions from $]0, 1]$ to the extended positive integers, i.e., the semigroup $\IZ^+ \cup \{+\infty\}$.
(Corresponding to a Hilbert module over $\C_0(]0, 1])$, i.e., a continuous field of Hilbert spaces over the
interval $]0, 1]$, one associates the function which at each point of $]0, 1]$ is equal to the dimension of the
fibre at that point.) Namely, one considers the particular functions $e_t=1_{]t, 1]},\ t\in [0, 1]$, and if
$\varphi$ and $\psi$ are two morphisms in $\cC u$ from $\cC u \C_0(]0, 1])$ to $\cC u A$, for each $r\ge 0$ one
says that the distance from $\varphi$ to $\psi$ is less than or equal to $r$ if $\varphi e_t \le\psi e_{t-r}$
and $\psi e_t\le \varphi e_{t-r}$ whenever $t\in [0, 1]$ and $t-r\ge 0$. (It will be convenient to include the
value $t=1$ even though $e_1 =0.$)

Let us verify that the number, say $d(\varphi, \psi)$, defined in this way satisfies the axioms for a metric. It
is clearly symmetric. It also clearly satisfies the triangle inequality: if $d(\varphi, \psi)\le r_1$, and
$d(\psi, \rho)\le r_2$, then, for any $t\in [0, 1]$ such that $t-(r_1 + r_2)\ge 0$, both
$$\psi e_t \le \psi e_{t-r_1}\le \rho e_{t-r_1-r_2}$$
and
$$\rho e_t\le \psi e_{t-r_2}\le\varphi e_{t-r_2-r_1};
$$
in other words, $d(\varphi, \rho)\le r_1 +r_2$. Finally, let us show that if $d(\varphi, \psi)=0$, i.e., if
$\varphi e_t=\psi e_t$ for every $t\in [0, 1]$, then $\varphi=\psi$. Our proof of this is rather indirect. (We
shall in particular use the hypothesis that $A$ has stable rank one.)

In fact, it was premature to try to introduce the metric to begin with---we must first prove that the semigroup
map under consideration is surjective!

Let, then, $\varphi$ be a homomorphism from $\cC u \C_0(]0, 1])$ to $\cC u A$, and let us show that $\varphi$
arises from a C*-algebra homomorphism from $\C_0(]0, 1])$ into $A\otimes \K$ (i.e., from an element of $\cT h
A)$.

The first step is to show that the decreasing family $(\varphi e_t)_{t\in [0, 1]}$ in $\cC u A$ may be realized
by a decreasing family $(E_t)_{t\in[0, 1]}$ of Hilbert $A$-modules---decreasing in the sense of inclusion, and
with the compact containment $\varphi e_s >> \varphi e_t$ in $\cC uA$ which holds whenever $s<t$ (because of the
relation $e_s>>e_t$ in $\cC u \C_0(]0, 1])$) realized by compact containment of $E_t$ as a subobject of $E_s$ in
the sense of [4], i.e., with the existence of a self-adjoint compact endomorphism of $E_s$ equal to the identity
on $E_t \subseteq E_s$.

Recall from  Theorem 1 of [4] that, in general, given two Hilbert C*-modules $E$ and $F$ over a C*-algebra $B$,
compact containment of the class [$E$] in the class [$F$], in the Cuntz semigroup $\cC u B$, written $[E]<<[F]$,
has two equivalent definitions, one involving only the order relation in $\cC u B$, and the other involving this
order relation and also the compact subobject containment relation recalled above---namely, it consists of the
requirement that there should exist a subobject $E' \subseteq F$ compactly contained as above, but with $[E']$
assumed not to be equal to $[E]$ but simply to majorize $[E]$ in $\cC u B$. In the case that $B$ has stable rank
one, by Theorem 3 of [4] the relation $[E]\le [E']$ in $\cC u B$ is equivalent to the relation that $E$ is
isomorphic to a subobject of $E'$, and so the relation $[E]<<[F]$ is equivalent just to the relation that $E$ is
isomorphic to a compactly contained subobject of $F$. (Note that the relation of compact containment for
subobjects, the definition of which was recalled above, is transitive, as compact endomorphisms of a subobject
extend canonically to a larger object, and in particular in a way preserving self-adjointness---see Theorem 6
below.)

To summarize, in the present case in which $A$ has stable rank one, if $E$ and $F$ are Hilbert $A$-modules such
that $[E]<< [F]$ in $\cC u A$, then $E$ is isomorphic to a compactly contained subobject of $F$. (It would be
convenient to know that the subobject of $F$ isomorphic to $E$ could be chosen to contain a given subobject of
$F$ which is isomorphic to a subobject of $E$. This might very well be true---but we don't  need it.)

Choose a Hilbert $A$-module $E_0$ such that $[E_0]=\varphi e_0$ in $\cC u A$. For each $n=0, 1, 2, \cdots$,
choose a decreasing sequence $E^n_0\supseteq E^n_{1\over 2^n}\supseteq\cdots \supseteq E^n_1=0$ of Hilbert
$A$-modules with $[E^n_t]=\varphi e_t$ in $\cC u A$ for $t={k\over 2^n},\ k=0, 1, \cdots, 2^n -1$, and with each
inclusion compact in the sense of [4], recalled above. (Make the choices one by one in decreasing order as
described in the preceding paragraph.) The case $n=0$ is the case considered at the beginning---let us write
$E^0_0=E_0$. In fact, in the general case one may suppose that $E^0_n=E_0$, since (by Theorem 3 of [4]) at least
$E^0_n$ is isomorphic to $E_0$---and we may replace the other modules in the decreasing sequence by their images
under this isomorphism---the compactness of the inclusions will be preserved.

For each fixed $n=0, 1, 2, \cdots$, and for each fixed $k=1, 2, \cdots, 2^n$, since $E^n_{k-1\over 2^n}$
compactly contains $E^n_{k\over 2^n}$ we may choose a self-adjoint compact  endomorphism of $E^n_{k-1\over
2^n}$, say $h_{n, k}$, which is the identity on $E^n_{k\over 2^n}$. Replacing $h_{n, k}$ by a continuous
function of it, we may suppose that it is positive and of norm at most one. Denote by $h_n$ the average of these
elements:
$$h_n:=2^{-n} \sum^{2^n}_{k=1} h_{n, k}.$$

\noindent Note that
$$0\le h_n - 2^{-n}\sum^{2^{n}}_{k=1} 1_{E_{k\over 2^n}}
\le 2^{-n},$$

\noindent where $1_{E_{k\over 2^n}}$ denotes the projection on $E_{k\over 2^n}$ in the bidual of the C*-algebra
of compact endomorphisms of $E_0$, and so any two choices of the sequence $h_{n, 1}, \cdots, h_{n, 2^n}$ in the
construction of $h_n$ yield the same result to within an error of $2^{-n}$.

Let us show that the sequence $(h_n)$ in the C*-algebra of compact endomorphisms of $E_0$ is Cauchy modulo inner
automorphisms (of the C*-algebra with unit adjoined).  For each $n$, note that the decreasing sequence $E_0 =
E_0^{n+1} \supseteq E^{n+1}_{{2 \over 2^{n+1}}} \supseteq \dots \supseteq E_1^{n+1} = 0$ of even-numbered terms
is exactly an alternative choice for the decreasing sequence $E_0 = E^n_0 \supseteq E^n_{{1 \over 2^n}}
\supseteq \dots \supseteq E^n_1 = 0$, and that the successive averages $(h_{n+1, 1} + h_{n+1, 2}) / 2, \cdots,
(h_{n+1, 2^{n+1} - 1} + h_{n+1, 2^{n+1}}) / 2$ constitute a particular choice for $h_{n, 1}, \cdots, h_{n,2^n}$
relative to this decreasing sequence of Hilbert $A$-modules, so that $h_{n+1}$, which is equal to the average of
the elements $(h_{n+1, 1} + h_{n+1, 2})/2, \cdots,$ is in particular just an alternative choice of $h_n$---with
respect to an alternative decreasing sequence of modules.  It is therefore sufficient, in order to prove that
the sequence $(h_n)$ is Cauchy modulo approximate unitary equivalence, to show that for each fixed $n$, elements
$h_n$ and $h'_n$ chosen with respect to two decreasing sequences $E_0 = E^n_0 \supseteq E^n_{{1 \over 2^n}}
\supseteq \dots \supseteq E^n_1 = 0$ and $E_0 = F_0 \supseteq F_{{1 \over 2^n}} \supseteq \dots \supseteq F_1 =
0$ of pairwise isomorphic Hilbert $A$-modules (with each one compactly contained in the preceding one) must be
approximately unitarily equivalent to within a tolerance which is summable with respect to $n$.  Note that it
was proved above that if the decreasing sequences are equal then $h'_n$ is approximately equal to $h_n$, to
within the (summable) tolerance $2^{-n}$.  Let us show that, when the decreasing sequences are not equal, the
elements $h_n$ and $h'_n$ are approximately unitarily equivalent (in the endomorphism C*-algebra of $E_0$) to
within the tolerance $4 \cdot 2^{-n} + \epsilon$ for any $\epsilon > 0$.

The enlargement of the tolerance (by over a factor of four) arises as follows.  We shall first consider a
particular choice for $h_{n-1}$, with respect to the even-numbered decreasing subsequence $E^n_0 \supseteq
E^n_{{2 \over 2^n}} \supseteq \dots \supseteq E^n_1 = 0$, which is then as shown above within $2^{-(n-1)} = 2
\cdot 2^{-n}$ of $h_n$.  We shall then show that this particular choice of $h_{n-1}$---which for clarity we
shall just call $g$---is, up to unitary equivalence in the C*-algebra of compact endomorphisms of $E^n_0 = E_0$
(with unit adjoined), arbitrarily close to an element $g'$ which is a particular choice for $h'_{n-1}$, i.e., an
element constructed as described earlier, with respect to the even-numbered decreasing sequence $E_0 = F_0
\supseteq F_{{2 \over 2^n}} \supseteq \dots \supseteq F_1 = 0$.  Note that $g'$ is not a fixed element, but
depends on the desired degree of approximation, and might better be denoted by $g'_{\epsilon}$.  In any case, as
shown above, $g'$ is within distance $2^{- (n-1)} = 2 \cdot 2^{-n}$ of $h'_n$, and so, with respect to a
suitable unitary $u$ depending on a given $\epsilon > 0$,
$$(*) \ \| h'_n - u h_n u^{-1} \| \le 2 \cdot 2^{-n} + \epsilon
+ 2 \cdot 2^{-n} = 4 \cdot 2^{-n} + \epsilon,$$ as desired.

In fact, we shall proceed to this inequality in a slightly different way from this, with an additional small
step which will be described clearly below.

By a particular choice of $h_{n-1}$, with respect to the decreasing sequence $E_0 = E^n_0 \supseteq E^n_{{2
\over 2^n}} \supseteq \dots \supseteq E^n_1 = 0$, we mean, of course, a particular choice of positive elements
$h_{n-1, k}$ of norm at most one, $k = 1, 2, \cdots, 2^{n-1}$, with $h_{n-1, k}$ a compact endomorphism of
$E^n_{{2 (k-1) \over 2^n}}$ for each $k = 1, 2, \cdots, 2^{n-1}$, equal to the identity on $E^n_{{2k \over
2^n}}$.  For each $k = 1, 2, \cdots, 2^{n-1}$, the element $h_{n, 2k-1}$ is such a choice; our choice of
$h_{n-1}$, then, will be the average of the elements $h_{n, 2 k - 1}, \ k = 1, 2, \cdots, 2^{n-1}$.  As shown
above, $\|h_n - h_{n-1} \| \le 2^{-(n-1)}$.

Analogously, by a particular choice of $h'_{n-1}$, with respect to the decreasing sequence $E_0 = F_0 \supseteq
F_{{2 \over 2^n}} \supseteq \dots \supseteq F_1 = 0$, we mean a particular choice of positive elements $h'_{n-1,
k}$ of norm at most one, $k = 1, 2, \cdots, 2^{n-1}$, with $h'_{{n-1}, k}$ a compact endomorphism of $F_{{2(k-1)
\over 2^n}}$ for each $k = 1, 2, \cdots, 2^{n-1}$, equal to the identity on $F_{{2 k \over 2^n}}$.  In other
words, a particular choice for $h'_{n-1}$ means the average of choices of such elements $h'_{n-1, 1}, h'_{n-1,
2}, \cdots, h'_{n-1, 2^{n-1}}$.  Recall that any such average is within distance $2^{-(n-1)}$ of $h'_n$.

It is at this point that we must introduce an additional step in the proof, as mentioned above.  Namely, it does
not seem possible to choose elements $h'_{n-1, 1}, \cdots, h'_{n-1, 2^{n-1}}$, approximately unitarily
equivalent to the elements $h_{n, 2 k-1}, \ k = 1, 2, \cdots, 2^{n-1}$, with exactly the properties stipulated
above---which would be the straightforward way to proceed; it seems to be possible only to do this with the
property that, for each $k = 1, 2, \cdots, 2^{n-1}$, the element $h'_{n-1, k}$ is equal to the identity on
$F_{{2k \over 2^n}}$ weakened to the property that $h'_{n-1, k}$ is just arbitrarily close to the identity on
$F_{{2k \over 2^n}}$.  This however is not serious; the average of the elements $h'_{n-1, 1}, h'_{n-1, 2},
\cdots, h'_{n-1, 2^{n-1}}$ is still arbitrarily close to being within distance $2^{-(n-1)}$ of $h'_n$.  In
particular, for given $\epsilon > 0$, for a suitable such choice, say $g$, one has

$$\|h'_n - g \| < 2^{-(n-1)} + \epsilon.$$

Note, furthermore, that to obtain this inequality we do not need that the elements $h'_{n-1, 1}, \cdots,
h'_{n-1, 2^{n-1}}$ belong exactly to the compact endomorphism algebras of $F_{{2 \over 2^k}}, \cdots, F_1$, but
only that, when they are considered as compact endomorphisms of $F_0$, they approximately belong to these
subalgebras.  (Note that we are taking Theorem 6 below for granted.) With this weakening of the requirements for
$h'_{n-1, 1}, \cdots, h'_{n-1, 2^{n-1}}$, as we shall show, they may be chosen to be exactly unitarily
equivalent to $h_{n,1}, h_{n,3}, \cdots, h_{n,2^n -1}$, with respect to a unitary element of the algebra of
compact endomorphisms of $E_0$ with unit adjoined, say $u$, and then the average, $g$, of these elements
(namely, $u h_{n,1} u^{-1}, u h_{n, 3} u^{-1}, \cdots, u h_{n, 2^n - 1} u^{-1}$) is exactly equal to the unitary
transform of the average, namely, of the particular choice of $h_{n-1}$ with respect to $E_0 = E^n_0 \supseteq
E^n_{{2 \over 2^n}} \supseteq \dots \supseteq E_1 = 0$ made above.  Hence, as $\|h_n - h_{n-1}\| \le
2^{-(n-1)}$,
$$\|g - u h_n u^{-1} \| \le 2^{-(n-1)} = 2 \cdot 2^{-n}.$$

As a result one obtains the inequality (*) unaltered.

To obtain a unitary $u$ with the properties specified above, i.e.~with $u h_{n, 2 k-1} u^{-1}$ approximately
contained in (the endomorphism algebra of) $F_{{2 (k-1) \over 2^n}}$ and approximately equal to the identity on
$F_{{2k \over 2^n}}$ for each $k=1, 2, \cdots, 2^{n-1}$, let us proceed step by step using Theorem 5 below.
Suppose, inductively, that for some $r = 1, 2, \cdots, 2^{n-1} - 1$, there exists an isomorphism $u_r$ from
$E_{{2(2^{n-1} - r) \over 2^n}}$ to $F_{{2(2^{n-1} - r) \over 2^n}}$ such that the transforms of the $r$
elements
$$h_{n, 2 (2^{n-1} - r) + 1}, h_{n, 2 (2^{n-1} - r) + 3}, \dots,
h_{n, 2 (2^{n-1} - r) + 2r-1} = h_{n, 2^n -1}$$ \noindent by $u_r$ are, respectively, approximately contained in
the compact endomorphism algebras of $F_{{2(2^{n-1} - r) \over 2^n}}, F_{{2(2^{n-1} - r) + 2 \over 2^n}},
\cdots, F_{{2^n - 2 \over 2^n}}$ (the first condition is trivial), and approximately equal to the identity on
$F_{{2(2^{n-1} - r) + 2 \over 2^n}}, \cdots, F_1,$ and let us find an  isomorphism $u_{r+1}$ of $E_{{2 (2^{n-1}
- (r+1)) \over 2^n}}$ onto $F_{{2 (2^{n-1} - (r+1)) \over 2^n}}$ with these same properties, with only an
arbitrarily small weakening of the approximations, and with the additional property that the transform of $h_{n,
2 (2^{n-1} - (r+1)) + 1}$ by $u_r$ is approximately contained in the compact endomorphism algebra of $F_{{2
(2^{n-1} - (r+1)) \over 2^n}}$ (this property is in fact trivial, but must be mentioned so that it is ensured to
hold with only a very small weakening at the next stage) and is approximately equal to the identity on $F_{{2
(2^{n-1} - (r+1)) + 2 \over 2^n}}$---to within an arbitrarily close degree of approximation. This is a simple
application of Theorem 5 with the extension $u_{r+1}$ of $u_r$ required to agree with $u_r$ approximately on
finite sets of module elements giving rise approximately to the compact endomorphisms $h_{n, 2 (2^{n-1} - r) +
1}, h_{n, 2 (2^{n-1} - r) + 3}, \cdots, h_{n, 2^{n-1}}$ and $u^{-1}_r h'_{n, 2 (2^{n-1} - (r+1)) + 2} u_r$ of
$E_{{2 (2^{(n-1)} - r) \over 2^n}}$. (Recall that by Theorem 6, below, compact endomorphisms of a subobject $M$
of an object $N$ are in a natural way also compact endomorphisms of $N$.) Since this last element is a compact
endomorphism of $E^n_{2 (2^{n-1} - (r+1))+1 \over 2^n}$, on which submodule $h_{n, 2 (2^{n-1} - (r+1)) +1}$ is
equal to the identity, so that $h_{n, 2 (2^{n-1}-(r+1)) + 1}$ acts as the identity on this last element, the
element $u_{r+1} h_{n, 2 (2^{n-1} - (r+1))+1} u^{-1}_{r+1}$ acts as the unit on $u_{r+1} u^{-1}_r h'_{n, 2
(2^{n-1} - (r+1)) + 2} u_r u^{-1}_{r+1}$ which is approximately equal to $h'_{n, 2 (2^{n-1} - (r+1)) +2}$ by
construction---which in turn is equal to the identity on $F_{{2 (2^{n-1} - (r+1)) + 2 \over 2^n}}$---and so the
transform of $h_{n, 2 (2^{n-1} - (r+1)) + 1}$ by $u_{r+1}$ is approximately equal to the identity on $F_{{2
(2^{n-1} - (r+1)) + 2 \over 2^n}}$, to within an arbitrarily close degree of approximation, as desired.  It
remains to note, first, that the existence of $u_1$ is trivial, since $F_1 = 0$ (and so any isomorphism from
$E^n_{{2 (2^{n-1} - 1) \over 2^n}}$ to $F_{{2 (2^{n-1} - 1) \over 2^n}}$ will do for $u_1$), and, second, that,
in the case $r = 2^{n-1} - 1, \ u_{r+1}$ is an isomorphism from $E^n_{{2 (2^{n-1} - 2^{n-1}) \over 2^n}} = E^n_0
= E_0$ to $F_{{2 (2^{n-1} - 2^{n-1}) \over 2^n}} = F_0 = E_0$, and so, again by Theorem 5 below, may be
approximated by an inner automorphism, not only on finitely many elements of $E_0$ but, as an automorphism of
the C*-algebra of compact endomorphisms of $E_0$, on the finite number of compact endomorphisms under
consideration.

This shows that the sequence $(h_n)$ is (summably) Cauchy modulo inner automorphisms. To complete the proof of
surjectivity, it remains to show that the limit, say $h$, in the C*-algebra of compact endomorphisms of $E_0$ of
a Cauchy sequence of unitary transforms $(u_n h_n u^{-1}_n)$ with $h_n$ constructed as above with respect to a
$2^n$-step decreasing sequence $E_0 = E^n_0 \supseteq E^n_{{1 \over 2^n}} \supseteq \dots \supseteq E^n_1 = 0$
with $[E^n_{{k \over 2^n}}] = \varphi e_{{k \over 2^n}}$ in $\cC u A$ corresponds to an element of $\cT h A$
giving rise to the given homomorphism $\varphi$ from $\cC u \text{C}_0 (]0,1])$ to $\cC u A$.  Note that a
positive element of norm at most one of the C*-algebra of compact endomorphisms of $E_0$ determines a
homomorphism of the C*-algebra $\text{C}_0 (]0,1])$ into this C*-algebra.

Let us show first that for each $t \in [0,1]$, the subobject of $E_0$ corresponding to $e_t(h)= 1_{]t,
1]}{(h)}$, i.e., the closed submodule generated by the images of the compact endomorphisms $f (h)$ with $f \in
\text{C}_0 (]t, 1])$, belongs to the class $\varphi e_t$ in $\cC u A$.

Fix $ t \in [0,1]$, and fix $\epsilon > 0$.  For $n$ sufficiently large, with $h_n \in \K (E_0)$ as above, up to
unitary equivalence in $\K (E_0)^{\smallsim}$ (the C*-algebra $\K (E_0)$ with unit adjoined), $h_n$ lies
strictly within distance $\epsilon$ of $h$.  By Corollary 9, below, for such $n$ (i.e., for $n$ such that, after
replacing $h_n$ by a unitarily equivalent element of $\K (E_0)$, one has $\|h - h_n \| < \epsilon$),
$$[e_{t+2 \epsilon} (h)] \le [e_{t + \epsilon} (h_n)] \le
[e_t (h)]$$ \noindent in $\cC u A$.  Furthermore, by Theorem 7, below, for sufficiently large $n$, and with the
restriction that both $t$ and $\epsilon$ are dyadic rational, as we shall show below,
$$\varphi e_{t + 2 \epsilon} \le [e_{t+\epsilon} (h_n)]
\le \varphi e_{t}$$ \noindent in $\cC u A$.

The application of Theorem 7 just alluded to is not quite as immediate as that of Corollary 9 which preceded
it---in particular the modules to apply it to must first be specified!  Given that both $t$ and $\epsilon$ are
dyadic rational, we may require that $n$ be large enough that both $t$ and $\epsilon$ are integral multiples of
$2^{-n}$.  Let us show that this is large enough.  By construction of $h_n$, we have Hilbert $A$-modules
$E^n_{t+2 \epsilon} \subseteq E^n_t \subseteq E_0$, representing $\varphi e_{t+2 \epsilon}$ and $\varphi e_t$
respectively, and, as $h_n = h^*_n \in \K (E_0)$, also $e_{t + \epsilon} (h_n)$ is a subobject of $E_0$; by
construction, $h_n \ge t + 2 \epsilon$ on $E^n_{t + 2 \epsilon}$ (i.e., $\langle(h_n - (t + 2 \epsilon)) \xi,
\xi\rangle_A \ \ge 0$ for any $\xi \in E^n_{t + 2 \epsilon}$), and, furthermore, $h_n \le t$ on any vector in
the bidual of $E''_0$ orthogonal to $E^n_t$ (i.e., $\langle(t - h_n) \xi, \xi\rangle_{A''} \ge 0$ for any such
vector $\xi$).  Of course, also $h_n \ge t + \epsilon$ on (the spectral subspace) $e_{t+ \epsilon} (h_n)$ and
$h_n \le t + \epsilon$ on any vector in $E''_0$ orthogonal to $e_{t + \epsilon} (h_n)$.  The desired
inequalities
$$\varphi e_{t + 2 \epsilon} = [E^n_{t + 2 \epsilon}] \le
[e_{t + \epsilon} (h_n) ]$$ \noindent and
$$[e_{t + \epsilon} (h_n) ] \le [E^n_t] = \varphi e_t$$
\noindent in $\cC u A$ now follow immediately from Theorem 7.

Combining the four inequalities established above (for $n$ sufficiently large), applied in the appropriate
order, for a given dyadic rational $\epsilon\!>\!0$ but with, in succession, $t + 2 \epsilon, \ t + 2 \epsilon,
\ t$, and $t$ in place of a given dyadic rational $t \in [0,1[$ (of course we may assume $t + 4 \epsilon \in
[0,1[$, but this is not necessary), we obtain the chain of inequalities
$$\varphi e_{t + 4 \epsilon} \le [e_{t + 3 \epsilon} (h_n)]
\le [e_{t + 2 \epsilon} (h)] \le [e_{t + \epsilon} (h_n)] \le \varphi e_t,$$ \noindent which on omitting the
terms involving $h_n$ becomes

$$\varphi e_{t + 4 \epsilon} \le [e_{t + 2 \epsilon} (h)]
\le \varphi e_t.$$

It is a little more convenient to view these two inequalities separately, for given dyadic $t$ and $\epsilon$,
as
$$\varphi e_{t + 2 \epsilon} \le [e_t (h)],$$
\noindent and
$$[e_{t + 2 \epsilon} (h) ] \le \varphi e_t. $$

>From the first, as $\epsilon$ tends to $0$, since $\varphi$ is a morphism in $\cC u$ and so preserves increasing
sequential suprema, we obtain
$$\varphi e_t \le [e_t (h)],$$
\noindent and from the second, similarly, as $[e_t (h)] = \psi e_t$ where $\psi$ denotes the homomorphism $\cC u
C_0 (]0,1]) \rightarrow \cC u A$ determined by $h$, we obtain
$$[e_t (h)] \le \varphi e_t.$$

In this way we obtain $\varphi e_t = [e_t (h)]$ at least for dyadic rational $t$, and it follows, again by
considering suprema over a decreasing sequence of values of $t$, that this equation holds, as desired, for all
values of $t$ in $[0,1[$.

It should perhaps be pointed out explicitly that the homomorphism from $\text{C}_0 (]0,1])$ to $\K (E_0)$
determined by $h$ does give rise to a homomorphism in the category $\cC u$ from $\cC u{C}_0 (]0,1])$ to $\cC u
A$---and not just, as it does in the first instance by functoriality, from $\cC u \text{C}_0 (]0,1])$ to $\cC u
\K (E_0)$.  As it happens, we do not need this at the moment, as the map $e_t \mapsto e_t (h) \subseteq E_0$,
defined directly above, clearly preserves increasing sequential suprema when these are calculated (on the
codomain side) in the lattice of subobjects of $E_0$, and it was shown in [4] (see proof of Theorem 1 of [4])
that the Cuntz class of the supremum of an increasing sequence of subobjects of a Hilbert C*-module is the
supremum of the classes of the subobjects---and this is all that is needed here.  (Recall that subobject means
countably generated closed submodule.) In fact, the map $e_t \mapsto [e_t (h)] \in \cC u A$ does factor through
a map $\cC u \K (E_0) \to \cC u A$, corresponding to the map of Hilbert modules taking the (right) Hilbert $\K
(E_0)$-module $H$ to the (right) Hilbert $A$-module
$$H \otimes_{\K (E_0)} E_0$$
\noindent (defined, as in a special case in [4]---cf.~also [15]---as the completion of the algebraic inner
product with respect to the norm derived from the natural $A$-valued inner product). Indeed, the functorially
defined element of $\cC u \K (E_0)$ corresponding to the module $e_t = \text{C}_0 (]t,1])$ over $\text{C}_0
(]0,1])$ is the class of the tensor product\break $e_t \otimes_{\text{C}_0 (]0,1])} \K (E_0)$, where $\K (E_0)$
is considered as a left Hilbert $\text{C}_0 (]0,1])$-module with respect to the element $h$, and the tensor
product of this with $E_0$,
$$e_t \otimes_{\C_0 (]0,1])} \K (E_0) \otimes_{\K (E_0)} E_0,$$
\noindent as a right A-module, is easily seen to be the subobject $e_t (h) \subseteq E_0$ defined earlier (the
first tensor product is the closed right ideal of $\K (E_0)$ generated by $\K (e_t (h))$, and the second is the
closure of the product of this with $E_0$, which is just the closure of $\K (e_t (h)) E_0$, which is just equal
to $e_t (h)$).  (Here we do not need that the stable rank of $A$ is one.)

The map $\cC u \K (E_0)\to \cC u A$ just described, arising from a natural map between modules (consisting of
tensoring  with the module $E_0$, considered as a two-sided Rieffel-Morita equivalence bimodule between $\K
(E_0)$ and $A$), may also be viewed as arising from a homomorphism from $\K(E_0)$ into $A\otimes \K$. Namely, by
Theorem 2 of [4], $E_0$ is isomorphic to a closed submodule of the Hilbert $A$-module direct sum of a countable
infinity of copies of the Hilbert $A$-module $A$, and the C*-algebra of compact endomorphisms of this is in a
natural way isomorphic to $A\otimes \K$. By Theorem 6 below, the algebra of compact endomorphisms of $E_0$ is
canonically then (in terms of a chosen embedding of $E_0$ in the infinite direct sum) embedded into $A\otimes
\K$, as a hereditary sub-C*-algebra. The corresponding map from $\cC u \K (E_0)$ to $\cC u A\otimes \K$ (which
is in a natural way isomorphic to $\cC u A$---see Appendix 6 of [4]) takes a Hilbert $\K(E_0)$-module $X$ into
the tensor product $X\otimes_{\K(E_0)} (A\otimes \K)$, a Hilbert $(A\otimes \K)$-module, and this is isomorphic
to $X\otimes_{\K(E_0)}\bigoplus^\infty_1  A$ and hence to $X\otimes_{\K(E_0)}\bigoplus^\infty_1 E_0$, and so
corresponds to the Hilbert $A$-module $X\otimes_{\K(E_0)}E_0$ in $\cC u A\ (\cong \cC u A\otimes\K)$, as
desired.

It follows that the map $e_t\mapsto [e_t(h)]\in \cC u A$ arises, by restriction to the set $\{e_t\}_{t\in[0, 1[}
\subseteq \cC u \C_0 (]0, 1])$, from a homomorphism $\C_0(]0, 1])\to A\otimes \K$, i.e., from an element, say
$\rho$, of the Thomsen semigroup of $A$. We wish to show that $\rho$ gives rise to $\varphi$, as in the
statement of the theorem. What we have shown so far is that $\rho$---or, rather, the homomorphism that $\rho$
gives rise to---and $\varphi$ agree on the subset $\{e_t\}_{t\in [0, 1[}$ of $\cC u \C_0(]0, 1])$. It is now
that we must assume that $\varphi$ extends to a map on $\Hom_1 (\cC u \C([0, 1]$, $\cC u (A\otimes \K)^{\sim})$,
as $\rho$ of course also does. The proof that, as a consequence, $\rho$ and $\varphi$ actually coincide is
surprisingly complex.

Let us show, then, that if $\varphi$ and $\psi$ are two homomorphisms in $\cC u$ from $\cC u \C([0, 1])$ to $\cC
u (A\otimes \K)^{\sim}$ which agree on 1 and on $e_t=1_{]t, 1]}$ for each $t\in  [0, 1]$ then they agree on all
of $\cC u \C([0, 1])$. By additivity and preservation of increasing sequential suprema, it is sufficient to show
that $\varphi$ and $\psi$ agree on the function $1_{]s, t[} \in \cC u \C([0, 1])$ for all $s, t \in \IR$. We
shall do this in two steps---first, the case $s<0$, and second the case $s\ge 0$.

By additivity and positivity (i.e., the property of preserving the order relation) of both $\varphi$ and $\psi$,
for any $t\in [0, 1]$ and  any $\epsilon >0$, \vglue -.1truein
\begin{eqnarray*}
\varphi(1_{[0, t-\epsilon[}) + \varphi(1_{]t-\epsilon, 1]})&=&
\varphi(1_{[0, t-\epsilon[}+1_{]t-\epsilon, 1]})\\
&\le& \varphi(1_{[0, 1]})=\psi(1_{[0, 1]})\\
&\le& \psi(1_{[0, t[}+1_{]t-\epsilon, 1]})\\
&=& \psi(1_{[0, t[})+\psi(1_{]t-\epsilon, 1]})\\
&=& \psi(1_{[0, t[})+\varphi(1_{]t-\epsilon, 1]}).
\end{eqnarray*}
More briefly,
$$\varphi(1_{[0, t-\epsilon[})+c\le d\le
\psi(1_{[0, t[})+c,$$ where $c=\varphi(1_{]t-\epsilon, 1]})$ and $d=\psi (1_{[0, 1]}).$ Since $1_{[0, 1]}$ is
compact in $\cC u \C([0, 1])$ (i.e., compactly contained in itself), by virtue of Theorem 10 below, and $\psi$
preserves compactness, $d$ is compact. Hence by Theorem 1 of [7], \vglue -.1truein
$$\varphi(1_{[0, t-\epsilon[})\le \psi
(1_{[0, t[}).$$ \vglue -.1truein
\smallskip

Hence as $\varphi$ preserves suprema of increasing sequences, \vglue -.1truein
$$\varphi(1_{[0, t[})\le \psi(1_{[0, t[}),$$\vglue -.1truein \smallskip \noindent
and so by symmetry (and antisymmetry!) $\varphi(1_{[0, t[})=\psi(1_{[0, t[})$, as desired.

Furthermore, as a consequence of this, we may proceed in a similar way for any $s\ge 0$ and $t\in [0, 1]$ with
$s<t$, and any $\epsilon >0$:
\begin{eqnarray*}
 \varphi(1_{[0, s+\epsilon[})+
  \varphi(1_{]s+\epsilon, t-\epsilon[}) +
  \varphi(1_{]t-\epsilon, 1]})
&=& \varphi(1_{[0, s+\epsilon[}+
  1_{]s+\epsilon, t-\epsilon[} +
  1_{]t-\epsilon, 1]})\\
&\le& \varphi(1_{[0, 1]})=
  \psi (1_{[0,1]})\\
&\le& \psi (1_{[0, s+\epsilon[}+
  1_{]s, t[} +
  1_{]t-\epsilon, 1]})\\
&=& \psi(1_{[0, s+\epsilon[})+
  \psi(1_{]s, t[})+
  \psi(1_{]t-\epsilon, 1]})\\
&=& \varphi(1_{[0, s+\epsilon[})+
  \psi (1_{]s, t[})+
  \varphi(1_{]t-\epsilon, 1]}).
\end{eqnarray*}
More briefly,
$$\varphi(1_{]s+\epsilon, t-\epsilon[})+c \le d\le\psi
(1_{]s, t[})+c,$$

\noindent where $c=\varphi(1_{[1, s+\epsilon[})+ \varphi(1_{]t-\epsilon, 1]})$ and $d=\varphi(1_{[0, 1]})$. As
before, as $d$ is compact (i.e., $d<<d$), by Theorem 1 of [7],
$$\varphi(1_{]s+\epsilon, t-\epsilon[})\le
\psi(1_{]s, t[}),$$ \noindent and hence as before
$$\varphi(1_{]s, t[})=\psi(1_{]s, t[}),$$
as desired.

It is now possible to show that the positive real-valued function
$$(\varphi, \psi)\mapsto d(\varphi, \psi)$$
defined above on pairs of homomorphisms from $\cC u \C_0(]0, 1])$ to $\cC u A$---assumed to extend to $\C([0,
1])$ as in the statement of the theorem!---is a metric. It remains to show that if $d(\varphi, \psi)=0,$ then
$\varphi=\psi$. By definition, $d(\varphi, \psi)=0$ just means that $\varphi(1_{]t, 1]})=\psi(1_{]t, 1]})$ for
all $t\in [0, 1]$. By what was just proved above, this implies $\varphi=\psi$.

Next, let us show that the map from $\cT h A$ to $\Hom (\cC u \C_0(]0, 1]), \cC u A)$ does not decrease
distances by more than the factor 8. Let $h$ and $k$ be positive elements of $A\otimes \K$, of norm at most one,
denote by $\varphi$ and $\psi$ the associated homomorphisms, and let us show that if $d(\varphi, \psi)\le r$
then $h$ and $k$ are approximately unitarily equivalent in $(A\otimes \K)^{\smallsim}$, to within the tolerance
$8r+\epsilon$ for any $\epsilon>0$.

Rather than considering $h$ and $k$ themselves, it seems to be necessary to deal with approximations to $h$ and
$k$, of the kind introduced in the construction of $h$ and $k$, above. Namely, for any $n=1, 2, \cdots$,
consider the open projections
$$e\ge e_0=1_{]0, 1]} (h)\ge e_1=1_{]{1\over n}, 1]} (h)\ge
\cdots, \ge e_{n-1}=1_{]{n-1\over n}, 1]} (h),$$ and, similarly, $f_l=1_{]{l\over n}, 1]} (k)$ for $0\le l\le
n-1$. With $g_l$ denoting the continuous function from $[0, 1]$ to $\IR$ which is zero from $0$ to ${l\over n}$,
one from $l+1\over n$ to $1$,
 and linear from $l\over n$ to ${l+1\over n}$, for
 each $l=0, 1, \cdots, n-1$, consider
 the finite sequences $g_l(h)$ and $g_l (k)$ in
 $A\otimes \K$, $l=1, 2, \cdots, n-1$, and
 note that (cf.~above) $h$ and $k$
 are within distance $1\over n$ of the
 averages of $g_0 (h), \cdots, g_{n-1} (h)$ and of
$g_0 (k), \cdots, g_{n-1} (k)$, respectively.

In fact, as earlier, it is necessary to note that $h$ and $k$ are within $1\over n$ of the average of $h_0, h_1,
\cdots, h_{n-1}$ and of $k_0, k_1, \cdots, k_{n-1}$, respectively, for any finite sequence $(h_l)$ and $(k_l)$
of positive elements of $A\otimes \K$ of norm at most one such that for each $l=0, 1, \cdots, n-1$, on the one
hand $e_l h_l=h_l$ and $f_l k_l =k_l$, and on the other hand $h_l e_{l+1}=e_{l+1}$ and $k_l f_{l+1}=f_{l+1}$
(where $e_n=f_n=0$). We shall in fact need the slightly more subtle statement (cf.~also earlier) that the
approximation of $k$ by the average of $k_0, k_1, \cdots, k_{n-1}$ still holds to within ${1\over n} + \epsilon$
for given $\epsilon>0$ if the equations $f_l k_l =k_l$ and $k_l f_{l+1}=f_{l+1}$ just hold approximately to
within a suitably small tolerance  depending on $\epsilon$ (for instance, $\epsilon\over 2$). Starting with a
given sequence $h_0, \cdots, h_{n-1}$ as above, for given $\delta>0$ we propose to find a unitary $u\in
(A\otimes \cK)\smallsim$ such that $uh_0 u^{-1},\cdots, uh_{n-1} u^{-1}$ satisfy the conditions just stated with
respect to $\delta$.

Actually, this proposal is fully appropriate only in the case $r=0$ (already considered earlier), and in the
case $r>0$ we must modify it, in a more or less straightforward way, as follows. Since $n$ is arbitrarily large,
it is sufficient to consider the case that $r={m\over n}$ for some $m=1, 2, \cdots, n-1$.

Choose an isomorphism $v_1$ from the right Hilbert module over $A\otimes \K$ determined by $f_{n-m}$, namely,
the largest closed right ideal on which $f_{n-m}$ acts on the left as the identity---let us denote this by
$F_{n-m}$---, to a subobject of the right Hilbert module $E_{n-2m}$ determined by $e_{n-2m}$. As above, the
existence of $v_1$ follows by Theorem 3 of [4] from the fact that (as $d(\varphi, \psi)\le r)$ the class of
$F_{n-m}$ in $\cC u A\otimes \K =\cC u A$ is majorized by that of $E_{n-2m}$. Since $f_{n-m} k_{n-m}=k_{n-m}$,
it follows that
$$e_{n-2m} v_1 k_{n-m}=v_1 k_{n-m},$$
and since $h_{n-3m} e_{n-2m}=e_{n-2m},$ we have
$$h_{n-3m} v_1 k_{n-m} =v_1 k_{n-m}.$$
Again by hypothesis (and by Theorem 3 of [4]), there exists an isomorphism $u_1$ from the (right) Hilbert module
$E_{n-3m}$ over $A\otimes \K$ determined by $e_{n-3m}$ to a subobject of the Hilbert module $F_{n-4m}$
determined by $f_{n-4m}$. By Theorem 5, below, $u_1$ may be chosen to extend $v_1^*$ approximately on finitely
many elements (note that $v^*_1$ is an isomorphism from a subobject of $E_{n-3m}$ to $F_{n-m}$, a subobject of
$F_{n-4m}$). In particular, considering $v_1 k_{n-m}$ as an element of $E_{n-3m}$ (on which $u_1$ just acts by
left multiplication), we may suppose that $u_1 v_1 k_{n-m}$ is arbitrarily close to $v_1^* v_1 k_{n-m}=k_{n-m}$,
and hence also (since $\| u_1\| \le 1$) that $u_1 h_{n-3m} u_1^* k_{n-m}$ is (arbitrarily) close to $k_{n-m}$.
(Note that $u_1 h_{n-3m} u_1^* k_{n-m}$ is close to $u_1 h_{n-3m} u^*_1(u_1 v_1 k_{n-m})$, which is equal to
$u_1 h_{n-3m} v_1 k_{n-m}$, and in fact to $u_1 v_1 k_{n-m}$ since $h_{n-3m} e_{n-2m}=e_{n-2m}$ and $e_{n-2m}
v_1 =v_1$, and is therefore close to $k_{n-m}$.)

Repeating this back-and-forth procedure as many times as possible (as many times as there is room for), starting
at the second stage with an isomorphism $v_2$ of $F_{n-5m}$ onto a subobject of $E_{n-6m}$, approximately
agreeing with $v_1$ on $k_{n-m}$, and then choosing an isomorphism $u_2$ from $E_{n-7m}$ to a subobject of
$F_{n-8m}$, approximately agreeing with $u_1$ on $v_1 k_{n-m}$ and also on $h_{n-3m}$ (these considered as
elements of $E_{n-3m}$, on which both $u_1$ and $u_2$ act by left multiplication)---and continuing in a similar
way at succeeding stages---, and at the very last stage jumping to an extension from all of $E_0$ to all of
$E_0$, where $E_0$ denotes the Hilbert $(A\otimes \K)$-module corresponding to $e_0$, we obtain an automorphism
$u$ of $E_0$ such that, for each $l=1, 2, \cdots$ for which this makes sense, i.e., for each $l\ge 1$ equal to 3
modulo 4, such that $n-lm\ge 0$, the element
$$uh_{n-lm} u^*$$
is arbitrarily small in norm on $1-f_{n-(l+1)m}$, and also arbitrarily close (in norm) to acting as the identity
on $k_{n-(l-2)m}$, and therefore also on $f_{n-(l-3)m}$ (as $k_{n-(l-2)m} f_{n-(l-3)m}=f_{n-(l-3)m}$, where
$f_n=0)\colon$
\[
\xymatrix{
        F_{n-m}\hspace{.5cm}\subseteq\hspace{-.7cm}&
        F_{n-4m}\hspace{.5cm}\subseteq\hspace{-.7cm}&
        F_{n-5m}\hspace{.5cm}\subseteq\hspace{-.7cm}&
        F_{n-8m}\hspace{.2cm}\subseteq\hspace{-.7cm}&
        \cdots\hspace{.2cm}\subseteq\hspace{-.7cm}&
        E_{0}\\
        E_{n-2m}\hspace{.5cm}\subseteq\ar@{<-}[u]_{v_{1}}\hspace{-.7cm}&
        E_{n-3m}\hspace{.5cm}\subseteq\ar@{->}[u]_{u_{1}}\hspace{-.7cm}&
        E_{n-6m}\hspace{.5cm}\subseteq\ar@{<-}[u]_{v_{2}}\hspace{-.7cm}&
        E_{n-7m}\hspace{.2cm}\subseteq\ar@{->}[u]_{u_{2}}\hspace{-.7cm}&
        \cdots\hspace{.2cm}\subseteq\hspace{-.7cm}&
        E_{0}\ar@{->}[u]_{u}
    }
\]
It follows that the average of the elements
$$uh_{n-lm}u^*,\ l=3, 7, 11, \cdots,$$
is arbitrarily close to being within distance ${4m\over n}=4r$ of $k$.

Note, furthermore, that as $h_{n-lm}$ is zero (in particular) on $1-e_{n-(l+1)m}$ and is the identity (in
particular) on $e_{n-(l-3)m}$, for each $l=3, 4, \cdots$ (where $e_n=0)$, the average of the elements
$$h_{n-lm},\ l=3, 7, 11, \cdots,$$
is within distance ${4m\over n}=4r$ of $h$. It follows that $uhu^*$ is arbitrarily close to being within
distance $8r$ of $k$, as desired.

This completes the proof that the map from $\cT h$ to $\Hom(\cC u \C_0(]0, 1]), \cC u A)$ under investigation
decreases distances by no more than a factor of $8$, with respect to the metrics under consideration. (In fact,
it may not decrease them strictly at all.) It remains to show that it does not strictly increase distances,
i.e., that it is a contraction. (It may in fact be an isometry.)

Suppose that $h$ and $k$ are positive elements of $A\otimes \K$ of norm at most one, and denote the
corresponding elements of $\Hom(\cC u \C(]0, 1]), \cC u A)$ by $\varphi$ and $\psi$. Suppose that $\|h-k\| \le
r$, and let us show that
$$d(\varphi, \psi)\le r,$$
i.e. that $\varphi e_t\le \psi e_{t-r}$ and $\varphi e_t\le \varphi e_{t-r}$ for all $t\in[0, 1]$ such that
$t-r\ge 0$, where $e_t=1_{]t, 1]}\in \cC u \C_0(]0, 1])$.

If $r\ge 1$ then the inequality $d(\varphi, \psi)\ge r$ holds vacuously (as does, incidentally, the inequality
 $\| h-k\| \le r)$.
If $r<1$, then from $\|h-k\| \le r$ it follows by Corollary 9, below, that $\varphi e_t\le \psi e_{t-r}$ and
$\psi e_t \le \varphi e_{t-r}$ for any $t\in [r, 1]$, i.e., that $d(\varphi, \psi)\le r$, as desired.
\end{proof}
\vspace{8pt}

{\bf 5.\ Theorem}\ (cf.~[4]). {\it{Let $A$ be a C*-algebra of stable rank one, and let $E_1, E_2$ and $F_1, F_2$
be pairs of isomorphic countably generated Hilbert $A$-modules with $E_1 \subseteq F_1$ and $E_2\subseteq F_2$.
It follows that for any choice of isomorphism $u$ from $E_1$ to $E_2$ there exists an isomorphism $v$ from $F_1$
to $F_2$ approximately extending $u$, i.e., agreeing arbitrarily closely with $u$ on each finite subset of $E_1$
(with respect to the norm in $F_2$).

In the case that $F_1$ and $F_2$ are equal, say to $F$, the approximate extension may be chosen to be inner,
i.e., to be an element of the C*-algebra of compact endomorphisms of $F$ with unit adjoined.

If the hypothesis that $F_2$ is isomorphic to $F_1$ is weakened to the hypothesis that $F_2$ contains a
subobject isomorphic to $F_1$, the conclusion still holds with $v$ just an isomorphism of $F_1$ onto some
subobject of $F_2$ (not necessarily a specified one).}}

\begin{proof}
Both the first and third statements follow from the second statement, by replacing $E_1$ and $F_1$ by their
image by an isomorphism of $F_1$ onto, respectively, $F_2$ or a subobject of $F_2$. (This reduces the first
statement immediately to the second statement. It reduces the third statement to the case that both $E_1$ and
$F_1$ are subobjects of $F_2$, in which case the extension of the given isomorphism $u:E_1\to E_2$ to an
isomorphism $v:F_1\to F'_1\subseteq F_2$ is effected by extending it (approximately) all the way to $F_2$!)

Consider, then the second statement. (This was dealt with earlier in [4], but only within the proof of Theorem 3
of [4]; let us repeat this proof for the convenience of the reader.) Let $u$ be an isomorphism between the
subobjects $E_1$ and $E_2$ of $F$. By  Proposition 1.3 of [17], $u$ may be approximated on a given finite subset
$S$ by a compact homomorphism of norm at most one from $E_1$ to $E_2$, say $x$---, which by Theorem 6 which
follows we may then consider as a compact endomorphism of $F$, with the same norm. Then in the C*-algebra $B$
obtained by adjoining a unit to the C*-algebra of compact endomorphisms of $F$, which is Rieffel-Morita
equivalent to a closed two-sided ideal of $A$ and therefore also of stable rank one, $x$ may be approximated in
norm by an invertible element, say $w$. We may suppose that also $w$ has norm at most one. Set $v=w(w^*
\omega)^{-{1\over2}}$, so that $v$ is a unitary element of $B$---and so, according to the terminology of [4], an
inner automorphism of $F$---, and let us show that $v$ approximately extends $u$, on the given finite subset $S$
of $E_1$. It is sufficient to show that $v$ approximately agrees with $w$ on any element of $F$ on which $w$ is
approximately isometric. Since $w$ has norm at most one, so also does $|w|=(w^* w)^{1\over 2}$. Hence, if
$\xi\in F$ is such that $\|w \xi\|$ is close to $\| \xi\|$, then, as shown in [4], also $|w| \xi$ is close to
$\xi$:
\begin{eqnarray*}
\langle(1-|w|)\xi, 1-|w| \xi\rangle &=&
   \langle \xi, (1-|w|)^2\xi\rangle\\
&\le& \langle \xi, (1-|w|^2)\xi\rangle=\langle \xi, (1-w^* w)\xi
\rangle\\
&=& \langle \xi, \xi\rangle -\langle w\xi, w\xi\rangle.
\end{eqnarray*}

Since $\|v\|=1$ (we exclude the case that $F=0$), from the fact that $\xi$ is close to $|w| \xi$ it follows that
$\xi$ is close to $v|w| \xi=w\xi$, as desired.
\end{proof}
\vspace{8pt}

{\bf 6.\ Theorem}\  (cf.~[3]).\
{\it{Let $A$ be a C*-algebra, and let $E$ and $F$ be right Hilbert $A$-modules. Let $G$ be a closed submodule of
$E$. The canonical extension to $E$ of a finite-rank homomorphism from $G$ to $F$ has the same norm. It follows
that there is a unique norm-preserving linear map from $\K (G, F)$ to $\K(E, F)$ compatible with this extension
process---let us refer to this as the canonical extension process for compact endomorphisms. In  particular, one
has a canonical injective *-homomorphism of C*-algebras $\K(G)\to \K(F)$. The image of this map is a hereditary
sub-C*-algebra.}}
\medskip

\begin{proof}
This follows immediately from the linking algebra description of a Hilbert C*-module given in [1] (see
Proposition 2.3 of [1])---applied to the Hilbert $A$-module $E$. The C*-algebra $\K(G)$ of compact endomorphisms
of $G$ is then revealed in a natural way as a hereditary sub-C*-algebra of the corresponding C*-algebra $\K(E)$
for $E$, and the elements of a approximate unit of norm one for $\K(G)$ therefore on the one hand map $E$ in a
norm-decreasing way into $G$, and on the other hand act approximately as the unit when composed with the
canonical extension of any finite-rank operator from $G$ to $F$---and so are approximately absorbed.
\end{proof}
\vspace{8pt}

{\bf 7.}\ \ One of the fundamental results concerning Murray-von Neumann comparability of projections in a von
Neumann algebra holds for countably generated Hilbert C*-modules.
\medskip

\noindent {\bf Theorem 7.1.}\ \ {\it{Let $A$ be a C*-algebra. Let $E$ be a (right) Hilbert $A$-module, and let
$E_1$ and $E_2$ be closed submodules of $E$. Assume that $E_1$ and $E_2$ are countably generated as Hilbert
$A$-modules. Suppose that, in the bidual $E''$ of $E$---a Hilbert C*-module over the bidual $A''$ of $A$---no
non-zero element of $E''_1$ is orthogonal to $E_2$ (equivalently, $E''_1 \cap E^{\bot}_2=0$). It follows that
$E_1$ is isomorphic (as a Hilbert $A$-module---i.e., in a way preserving the $A$-valued inner product) to a
closed submodule of $E_2$.}}

\begin{proof}
Choose strictly positive elements $h$ and $k$ in the C*-algebras of compact endomorphisms of $E_1$ and $E_2$
respectively (which exist as $E_1$ and $E_2$ are countably generated by Corollary 1.1.25  of [15]), and let us
show that $hk^2h$ is also a strictly positive element of $\K(E_1)$---which it belongs to as this is a hereditary
sub-C*-algebra of $\K(E)$ by Theorem 6.

Let $f$ be a positive functional on $\K(E_1)$ which is zero on $hk^2 h$, and let us show that $f$ is equal to
zero. The hypothesis $E''_1\cap E^{\bot}_2 =0$ may be interpreted in view of the linking algebra description of
$E$ given in Proposition 2.3 of [1] as the condition that any positive functional on $\K (E)$ which is zero on
the hereditary sub-C*-algebra $\K(E_2)$ is also zero on $\K(E_1)$. The positive functional $\K (E)\ni x\mapsto
f(hxh)$ is zero on $k^2$ which is a strictly positive element of $\K(E_2)$, and it is therefore zero on all of
$\K(E_2)$. Hence by hypothesis it is zero on $\K (E_1)$, and in particular, $f(h^3)=0$. Since $h^3$ is a
strictly positive element of $\K (E_1)$, it follows that $f=0$.

This shows that $hk^2h$ is a strictly positive element of $\K (E_1)$. In other words, the closure of $hkE$ is
equal to $E_1$. It follows that the partially isometric part of the compact endomorphism $kh$ of $E$---which is
in any case an isomorphism from the closure of the range of its adjoint to the closure of its range---and so in
the present circumstances an isomorphism from the closure of $hkE$ to the closure of $khE$---is an isomorphism
from $E_1$ to a subobject of $E_2$, as desired.
\end{proof}
\vspace{8pt}

{\bf 8.}\ \ The following consequence of Theorem 7 is also a consequence of [13]; we include the present more
elementary proof as the result is used in the proof of Theorem 4, above. (More precisely, it is used in the
proof of Theorem 1 of [10], which in turn is used in the proof of Theorem 4 above.)
\smallskip

\noindent {\bf Corollary 8.1.}\ \ {\it{Let $A$ be a C*-algebra, and let $E$ and $F$ be Hilbert C*-modules over
$A$. Let $k:E\to F$ be a bounded $A$-linear map from $E$ to $F$. Let $E_0$ be a countably generated closed
submodule of $E$, and suppose that $k$ is bounded below on $E_0$. It follows that $E_0$ is isomorphic as a
Hilbert $A$-module to a closed submodule of the closed submodule $kE_0$ of $F$.}}

\begin{proof}
Choose a strictly positive element $h$ of $\K(E_0)$ (see Corollary 1.1.25 of [15]). In the bidual $E''$, suppose
that $\xi \in E''_0$ is orthogonal to the range of $(kh)^*$(note that $kh$ is compact, hence adjointable). In
other words, $\langle \xi, (kh)^*\eta\rangle_{A''}=0$ for all $\eta$ in (the bidual of) $F$, with respect to the
$A''$-valued inner product in the bidual of $E$, considered as a Hilbert C*-module over the bidual $A''$ of $A$.
Then $k h\xi=0$, from which it follows as $k$ is bounded below on $E_0$, and therefore also on $E''_0$, and that
$h$ is injective on $E''_0$, that $\xi=0$.

It follows by Theorem 7, above,  that $E_0$ is isomorphic as a Hilbert C*-module to a closed submodule of the
closure of the image of $(kh)^*$ (note that this closed submodule is countably generated by Corollary 1.1.25 of
[15], as $(kh)^*$ is compact). Since, for any compact (or even just adjointable) homomorphism $r$ between
Hilbert C*-modules, the closures of the images of $r$ and of $r^*$ are isomorphic as Hilbert C*-modules (by
means of the partially isometric part of $r$), it follows that $E_0$ is isomorphic to a closed submodule of the
closure of $khE$, i.e., since $hE$ is dense in $E_0$ and $k$ is bounded below on $E_0$, to a closed submodule of
$kE_0$.

(Note that, by Theorem 4.1 of [13], $E_0$ is in fact isomorphic to  $k E_0$ itself, not just to a closed
submodule.)
\end{proof}
\vspace{8pt}

{\bf 9.}\ \ The following consequence of Theorem 7 is also used in the proof of Theorem 4, above.
\smallskip

\noindent {\bf Corollary 9.1.}\ \ {\it{Let $h_1$ and $h_2$ be self-adjoint compact endomorphisms of the Hilbert
C*-module $E$ over $A$, and suppose that $\|h_1-h_2\| \le \varepsilon$. It follows that, for every $t\in\IR$,
the closed submodule of $E$ corresponding to the spectral projection $1_{]t, +\infty]} (h_1)$ for $h_1$ (i.e.,
the range of this open projection) is isomorphic, as a Hilbert $A$-module, to a closed submodule of the closed
submodule of $E$ corresponding to the spectral projection $1_{]t-\varepsilon, +\infty]} (h_2)$ for $h_2$.}}

\begin{proof}
With $E_1$ and $E_2$ denoting the two closed submodules of $E$ referred to in the statement of the
theorem---which are countably generated by Corollary 1.1.25 of [15] as their compact endomorphism C*-algebras
have countable approximate units---note that for fixed $t\in \IR$,
$$\langle \xi, h_1 \xi\rangle_{A''} >
t\langle \xi, \xi\rangle_{A''}$$ for every non-zero $\xi\in E''_1$, and
$$\langle\xi, h_1 \xi\rangle_{A''}\le
\langle \xi, (h_2+\varepsilon) \xi\rangle_{A''}\le t\langle \xi, \xi\rangle_{A''}$$ for every $\xi \in
E^\bot_2$, from which follows $E''_1\cap E^\bot_2=0$. Hence by Theorem 7, $E_1$ is isomorphic to a closed
submodule of $E_2$.
\end{proof}
\vspace{8pt}

{\bf 10.}\ \ The following computation of $\cC u \C([0, 1])$ and of $\cC u \C_0 (]0, 1])$ is fundamental to
Theorems 4 and 11.
\medskip

\noindent {\bf Theorem 10.1.}\ ([14]).\ \ {\it{$\cC u \C_0(]0, 1])$ is isomorphic to the ordered semigroup of
lower semicontinuous extended positive integer valued functions on the half-open  interval $]0, 1]$, by means of
the map which to any Hilbert C*-module over $\C_0 (]0, 1])$ associates its dimension function, i.e., the
function on the spectrum of $\C_0(]0, 1])$ the value of which at any pure state is the dimension of the Hilbert
space constructed from the given Hilbert C*-module by composing the pure state with the $\C_0(]0, 1])$-valued
inner product.}}

\begin{proof}
This follows from [14] and [3]. The case of finitely generated Hilbert C*-modules could also be dealt with using
the result that the eigenvalue pattern of a self-adjoint element of M$_n (\C(]0, 1]))$ (or of $\rM_n(\C_0(]0,
1]))$) determines its approximate unitary equivalence class, together with Theorem 3 and Appendix 6 of [4] which
show that our present definition of the Cuntz semigroup in the stable rank one case (isomorphism classes of
Hilbert C*-modules!)~is equivalent (at least as far as finitely generated Hilbert C*-modules are concerned) with
that of [5].
\end{proof}
\vspace{8pt}

{\bf 11.}\ \ If follows from  Theorem 4 that the Thomsen semigroup in fact determines the Cuntz semigroup.
\smallskip

\noindent {\bf Theorem 11.1.}\ \ {\it{If $A$ has stable rank one, then the ordered semigroup $\cC u A$ is
determined by the semigroup $\Hom (\cC u \C_0(]0, 1]), \cC u A)$, as the set of all images under homomorphisms
of the constant function $1_{]0, 1]}\in \cC u \C_0(]0, 1])$, with the addition inherited from addition of
homomorphisms, and the order relation inherited from the inequality $1_{]{1\over 2}, 1]} << 1_{]0, 1]}$ in $\cC
u \C(]0, 1])$, in the sense that $x<< y$ in $\cC u A$ (recall that the order relation in $\cC u A$ is determined
by the relation $<<$) if and only if $x$ and $y$ are simultaneously the images  of $1_{]{1\over 2}, 1]}$ and
$1_{]0, 1]}$ by a homomorphism $\cC u \C_0(]0, 1])\to \cC u A$ in the category $\cC u$. It is sufficient to
consider elements of $\Hom (\cC u \C_0(]0, 1]), \cC u A)$ which extend to elements of $\Hom_1 (\cC u C([0, 1]),
\cC u (A\otimes \K)^{\smallsim})$ (i.e.~by Theorem 4, which arise from elements of $\cT h A$).}}

\begin{proof}
It is sufficient to prove that if $x$ and $y$ are elements of $\cC u A$ with $x<<y$ then there is a homomorphism
$\varphi: \cC u \C_0(]0, 1])\to \cC u A$ taking $1_{]{1\over 2}, 1]}$ into $x$ and $1_{]0, 1]}$ into
$y$---arising from an element of $\cT h A$.

Since, whenever $x<<y$ in an ordered semigroup in the category $\cC u$, there exists a rapidly increasing
sequence $x_1<<z_1<<z_2<<\cdots <<y$ with $\sup z_n=y$, we may choose a family $(z_t)_{t\in \IZ [{1\over 2}]\cap
[0, 1]}$ in $\cC u A$ such that $z_1=0$, $z_{1\over 2}=x$, and $z_0=y$, such that $z_s << z_t$ whenever $s>t$,
and such that $\sup z_{t_n}=z_t$ whenever $t_1\ge t_2 \ge\cdots$ is a sequence in $\IZ[{1\over 2}]\cap [0, 1]$
decreasing to $t\in \IZ [{1\over 2}]$. Then we may extend this to a family $(z_t)_{t\in [0, 1]}$ in $\cC u A$
such that for every decreasing sequence $t_1\ge t_2\ge \cdots$ in $[0, 1],\ z_{\inf t_n}=\sup z_{t_n}$ in $\cC u
A$. As shown in the proof of Theorem 4, the family $(z_t)$ arises from a positive compact endomorphism of a
certain countably generated Hilbert C*-module over $A$, which may be taken to be that (unique up to isomorphism
as $A$ has stable rank one) corresponding to $y$, and hence from a C*-algebra homomorphism $\C_0 (]0, 1])\to
A\otimes \K$, taking $1_{]{1\over 2}, 1]}$ into $x$ and $1_{]0, 1]}$ into $y$ as desired.

In fact, the construction of the family $(z_t)$ as above can clearly be carried out in any object in the
category $\cC u$ (beginning with a pair of elements $x<<y$), in place of the particular object $\cC u A$ under
consideration. It would be interesting to know if, for any such object $Z$, the family $(z_t)$ can be extended
to a homomorphism
$$\cC u \C_0 (]0, 1]) \to Z,$$
taking $1_{]t, 1]}$ into $z_t$ for each $t\in [0, 1]$.
\end{proof}
\vspace{8pt}

{\bf 12.}\ \ The classification of (unital) AI algebras in terms of the Cuntz semigroup (together with the class
of the unit) that follows (virtually) immediately from [19] and Theorem 4 above can be extended slightly as
follows. Note that by Corollary 4 of [3] the class of C*-algebras considered in the following theorem is closed
under stable isomorphism.
\smallskip

\noindent {\bf Theorem 12.1.}\ \ {\it{Let $A$ and $B$ be the C*-algebra inductive limits of sequences of
separable C*-algebras of continuous trace with spectrum a finite union of closed or half-open intervals in the
real line. Let $\varphi_{\cC u}:\ \cC u A\to \cC u B$ be an isomorphism of ordered semigroups taking the class
in $\cC u A$ of the (countably generated) Hilbert $A$-module $A$ into the class in $\cC u B$ of the Hilbert
$B$-module $B$.

It follows that there exists an isomorphism $\varphi:A\to B$ giving rise to the given isomorphism $\varphi_{\cC
u}$.}}

\begin{proof}
By [14], up to stable isomorphism the class of C*-algebras under consideration is the class of inductive limits
of sequences of the tensor product of $\K$ with a finite direct sum of copies of $\C([0, 1])$ and $\C_0(]0,
1])$. Let us establish the conclusion of the theorem first in the case that $A$ and $B$ belong to this class
(i.e., are stable).

By Theorems 4 and 11 the isomorphism $\varphi_{\cC u}$ gives rise to an isomorphism $\varphi_{\cT h}: \cT h A\to
\cT h B$, which is not a priori isometric with respect to the natural metric (although it is an immediate
consequence of the conclusion of the theorem that it must be!), but which, together with its inverse, is at
least uniformly continuous. This is enough to carry through the proof of Theorem 1.5 of [19]---with minor
modifications to adapt it to the present setting, with mixed building blocks (spectrum  $[0, 1]$ or $]0, 1]$),
assumed to be stable rather than unital (in fact, they cannot be unital when the spectrum is noncompact), and
with, a priori, only the uniform structure underlying the metric on the Thomsen semigroup known to be preserved
by the isomorphism. Let us in fact carry out this proof (since the present setting does differ, in two ways,
from that of [19]).

Let $A_1\to A_2 \to\cdots$ and $B_1\to B_2 \to \cdots $ be sequences of finite direct sums of copies of the
tensor product of $\C([0, 1])$ or $\C_0 (]0, 1])$ with $\K$, with inductive limits $A$ and $B$ respectively. The
proof is an intertwining argument, analogous to that of [5] and [6], but with the difference that, as in [7] and
[19], the intertwining already at the level of the invariant (not just at the level of the algebras) is only
approximate---that is, the resulting diagram is only approximately commutative. As in [7] and [19], the
intertwining at the level of the invariant is by means of maps that are metric space contractions, even though
the given isomorphism between the invariants of $A$ and $B$ is not known to have this property, but just to
preserve the uniform structure! (This would appear to be crucial in order to achieve stable control of the
discrepancies in the intertwining.)

The basic step in the intertwining is to consider the finite stage algebra $A_1$ (for instance) and to compose
the canonical map $\cT h A_1 \to \cT h A$ with the given map $\cT h A\to \cT h B$, to obtain a uniformly
continuous map $\cT h A_1 \to \cT h B$. Since it preserves the module structure of $\cT h$ (over the
endomorphism monoid of $\C_0 (]0, 1])$), and both $A_1$ and $B$ are stable, it follows that this map in fact
arises from a C*-algebra homomorphism $A_1 \to B$. In particular, not only is it a contraction, but also this
holds without using the hypothesis of uniform continuity.

It follows at this stage from uniform continuity of the map $\cT h A\to \cT h B$, owing to the fact (depending
on the stability of  relations in
 $\C_0 (]0, 1])$---cf. [19])
that $\cT h A$ is the inductive limit of the sequence $\cT h A_1 \to \cT h A_2 \to \cdots,$ in the category of
complete metric spaces (also with additional structure, which is carried automatically by continuity to the
inductive limit), that the map $\cT h A\to \cT h B$ is a contraction, and hence by symmetry an isometry.

In the case that all the minimal direct summands of the finite-stage algebras $A_i$ and $B_i$ have spectrum the
closed interval, the conclusion of the theorem may now be obtained just as in [19], the fact that all the
C*-algebras being considered are now stable instead of unital not affecting the proof of Theorem 1.5 of  [19] in
any material way---all three stages of the proof---intertwining at the level of the invariant, lifting the
individual maps of this intertwining to C*-algebra maps between building blocks, and modification of these maps
to obtain an approximately commutative diagram in the sense of Theorems 2.2 and 2.1 of [6]---being still valid
without remark.

In the case that some of the closed intervals are replaced by half-open intervals, the same steps are still
valid, but the following remark must be made. Namely, while the first two of the three steps of the intertwining
are still seen to be valid in the same way as in [19], the question of approximate unitary equivalence of two
C*-algebra maps from $A_1$ to, say, $A_2$, which agree approximately at the level of $\cT h$, requires some
comment. It is enough to consider the case that $A_1$ is equal to $\C_0 (]0, 1])\otimes \K$. It is in fact
enough to consider the case that $A_1$ is equal to $\C_0(]0, 1])\otimes \rM_n$ for some $n$ (in other words, to
pass to such a subalgebra). In the case $n=1$, approximate unitary equivalence (to within a given tolerance)
follows immediately from closeness at the level of the invariant (to within the same tolerance). Already in the
case $n=2$, a new difficulty arises, owing to the fact that the algebra $\C_0(]0, 1])$ is non-unital, and so
$\rM_n$ is not a subalgebra of $A_1$. Nevertheless, if $e_{12}$ denotes the (upper right hand) off-diagonal
generator for $\rM_2$, then the element $h\otimes e_{12}$ of $A_1=\C_0(]0, 1])\otimes \rM_2$, where $h$ denotes
the canonical generator $t\mapsto t$ of $\C_0(]0, 1])$, which generates $A_1$, is determined up to approximate
unitary equivalence, to within a given tolerance, in $A_2$, or in any C*-algebra $D$ of stable rank one into
which $A_1$ is mapped, by knowing the image of $h\otimes e_{11}$---or even just by knowing this approximately,
to within a sufficiently small tolerance---which depends only on the given tolerance, and not on $D$ or the maps
$A_1\to D$. This resolves the difficulty, but the proof of this assertion is a slightly subtle application of
the stable rank one isomorphism extension result for (countably generated) Hilbert C*-modules---Theorem 5
above---, which we must now examine. It will be clear from our examination of this case, i.e. the case $n=2$,
that a similar assertion holds also in the case of any finite $n$ (concerning the elements $h\otimes e_{1i}$
with $i=2, 3, \cdots, n$).

Replacing $D$ by the asymptotic sequence algebra with respect to a sequence of C*-algebras $D_1, D_2, \cdots$ of
stable rank one, i.e., the C*-algebra of bounded sequences ($d_i$) with $d_i \in D_i$ modulo the closed
two-sided ideal of null sequences, we may reduce the question involving only approximate knowledge of the image
of $h\otimes e_{11}$ to the question in the case the image of $h\otimes e_{11}$ in $D$ is known exactly. Thus,
one has two elements $x_1$ and $x_2$ of norm at most one and of square zero such that the elements $x^*_1  x_1$
and $x^*_2  x_2$ are equal. There is then a Hilbert $D$-module isomorphism of the closed right ideal of $D$
generated by $x_1$ and $x^*_1$ onto the closed right ideal of $D$ generated by $x_2$ and $x^*_2$, taking the
compact endomorphism $x_1$ onto the compact endomorphism $x_2$, and hence by Theorem 5 there is a unitary
element $u$ of the C*-algebra obtained by adjoining a unit to the C*-algebra of compact endomorphisms of the
(right) Hilbert $D$-module $D$---which is just $D$!---which gives rise approximately to this isomorphism between
closed submodules, and in particular takes the compact endomorphism $x_1$ approximately onto the compact
endomorphism $x_2$, in other words, such that, when $x_1$ and $x_2$ are considered as elements of $D$
(cf.~Theorem 6 above), the elements $u x_1 u^*$ and $x_2$ of $D$ are arbitrarily close, as desired.

Finally, consider the general case, in which $A$ and $B$ are no longer assumed to be stable. Considering
$A\otimes \K$ and $B\otimes \K$, which by Appendix 6 of [3] have, in a natural way, the same Cuntz semigroups as
$A$ and $B$, we have by what was shown above an isomorphism of $A\otimes \K$ with $B\otimes \K$ giving rise to
the given isomorphism between $\cC u A$ and $\cC u B$. The hereditary sub C*-algebra $A\otimes e$ of $A\otimes
\K$, where $e$ is a minimal non-zero projection in $\K$, is then taken onto a full hereditary sub-C*-algebra of
$B\otimes \K$, and the Hilbert module over $B\otimes \K$ corresponding to it, i.e., the closed two-sided ideal
of $B\otimes \K$ generated by it, is isomorphic to that corresponding to $B\otimes e$, i.e., to the closed right
ideal generated by $B\otimes e$. (Recall that this isomorphism of Hilbert $B\otimes \K$-modules, instead of just
equality of the Cuntz equivalence classes, is owing to the fact that $A$ has stable rank one.) This amounts to
an isomorphism of the hereditary sub-C*-algebras in question, giving rise to the identity map on the Cuntz
semigroup, identified in each case with that of $B\otimes \K$. Combining this with the isomorphism of $A\otimes
e$ onto the subalgebra of $B\otimes \K$, one has an isomorphism from $A=A\otimes e$ to $B=B\otimes e$ giving
rise to the given isomorphism from $\cC u A$ to $\cC u B$. (Note that this reduction of the general case to the
stable case is valid for arbitrary separable C*-algebras of stable rank one---separable so that the Hilbert
C*-modules $A$ and $B$ are countably generated.)
\end{proof}

\textbf{Added August 25, 2007 } Recently, R\o rdam and Winter (in \emph{The Jiang-Su algebra revisited}, in
preparation) have obtained a strengthening of the cancellation result for the Cuntz semigroup of a C*-algebra of
stable rank one established in [7], and used in the proof of Theorem 4 of the present paper. This allows us to
simplify the statements of Theorems 4 and 11. More specifically, it allows us to show that the two semigroups of
Cuntz semigroup homomorphisms discussed in Section 4---and appearing in the statements of these two
theorems---associated, on the one hand, to the C*-algebra $\C_0(]0, 1])$ and a given stable C*-algebra $A$ of
stable rank one and, on the other hand, to the unitalizations of these two C*-algebras---the homomorphisms in
this case required to preserve the units---, are exactly the same. (Accordingly, the two theorems may be stated
without mentioning the algebras with units adjoined.)

To show that the two semigroups are the same, in the natural sense---that the homomorphisms between the Cuntz
semigroups of the two algebras before the adjunction of units are in natural bijective correspondence with the
unit-preserving homomorphisms after the adjunction of units, the correspondence between two homomorphisms just
consisting of their compatibility in the obvious sense---it is most convenient to prove, in addition to Theorem
4, the modified statement of Theorem 4, in which the non-unital algebras appear instead of the unitalized ones.

Let us show, then, that if $\varphi$ and $\psi$ are two homomorphisms in $\cC u$ from $\cC u \C([0, 1])$ to $\cC
u A$ which agree on $e_t=1_{]t, 1]}$ for each $t\in  [0, 1]$ then they agree on all of $\cC u \C([0, 1])$. By
additivity and preservation of increasing sequential suprema, it is sufficient to show that $\varphi$ and $\psi$
agree on the function $1_{]s, t[} \in \cC u \C([0, 1])$ for all $s, t \in \R$. By additivity and positivity
(i.e., the property of preserving the order relation) of both $\varphi$ and $\psi$, for any $s, t \in \R$ with
$s < t$, and any $\epsilon > 0$, \vglue -.1truein
\begin{eqnarray*}
\varphi(1_{]s+\epsilon, t-\epsilon[}) + \varphi(1_{]t-\epsilon,1]}) &=&
\varphi(1_{]s+\epsilon, t-\epsilon[} + 1_{]t-\epsilon, 1]})\\
&\le&\varphi(1_{]s+\epsilon, 1]})\\
&<<&\varphi(1_{]s, 1]}) = \psi(1_{]s, 1]})\\
&\le& \psi(1_{]s, t[}+1_{]t-\epsilon, 1]})\\
&=& \psi(1_{]s, t[})+\psi(1_{]t-\epsilon, 1]})\\
&=& \psi(1_{]s, t[})+\varphi(1_{]t-\epsilon, 1]}),
\end{eqnarray*}
since $1_{]s+\epsilon]} << 1_{]s, 1]}$ in $\cC u \C([0, 1])$. More briefly,
$$\varphi(1_{]s+\epsilon, t-\epsilon]}) + c << \psi(1_{]s, t[}) + c,$$
where $c = \varphi(1_{]t-\epsilon, 1]})$. Hence, by Proposition 4.3 of [18],
$$\varphi(1_{]s+\epsilon, t-\epsilon[}) \le \psi(1_{]s, t[}).$$

Hence as $\varphi$ preserves suprema of increasing sequences,
$$\varphi(1_{]s, t[}) \le \psi(1_{]s, t[}),$$
and also by symmetry (and antisymmetry) $\varphi(1_{]s, t[}) = \psi(1_{]s, t[})$, as desired.

Here, in quoting Proposition 4.3 of [18], we actually mean not the literal statement of this result, but the
following implication (which involves just three elements $a, b, c$ of the Cuntz semigroup)$\colon$
$$a + c << b + c \Rightarrow a << b$$

(This is easily seen by the results of [4]---without the hypothesis of stable rank one---to be equivalent to the
following implication$\colon$
$$a + c \le b + c', c' << c \Rightarrow a << b$$
Again because of the results of [4] (in the general case) this is readily seen to be equivalent to the
implication of Proposition 4.3 of [18], in which $c'$ is $[(x - \epsilon)_{+}]$ for some $x \geq 0$ with $[x]
= c.$)



\bigskip Department of Mathematics,\\University of Toronto,\\Toronto, Canada, M5S 2E4
\smallskip\\ciuperca@math.toronto.edu\\elliott@math.toronto.edu

\end{document}